\documentclass{amsproc}
\usepackage{amsmath, amssymb, amsthm, fullpage}
\usepackage[all]{xypic}
\newtheorem{thm}[equation]{Theorem}
\newtheorem{prop}[equation]{Proposition}

\theoremstyle{definition}
\newtheorem{defn}[equation]{Definition}

\theoremstyle{remark}
\newtheorem{exam}[equation]{Example}

\newtheorem{ntn}[equation]{Notation}

\newtheorem{rem}[equation]{Remark}

\makeatletter
\renewcommand{\subsection}{\@startsection{subsection}{2}{0pt}{-3ex
plus -1ex minus -0.2ex}{-2mm plus -0pt minus
-2pt}{\normalfont\bfseries}} \makeatother

\numberwithin{equation}{section}

\newdir^{ (}{{}*!/-5pt/@^{(}}

\newcommand{\beq}{\begin{equation}\label}
\newcommand{\eeq}{\end{equation}}

\newcommand{\sgn}{\operatorname{sign}}
\newcommand{\onto}{\twoheadrightarrow}

\newcommand{\Hom}{\operatorname{Hom}}

\newcommand{\Ind}{\operatorname{Ind}}
\newcommand{\SSS}{\mathbb{S}}

\newcommand{\J}{\mathcal{J}}

\newcommand{\cO}{\mathcal O}
\newcommand{\kk}{\mathbf{k}}
\newcommand{\ot}{\otimes}
\newcommand{\End}{\operatorname{End}}
\newcommand{\dq}{\overline{Q}}
\newcommand{\Id}{\mathrm{Id}}
\newcommand{\id}{\mathrm{id}}
\newcommand{\F}{\mathcal{F}}
\newcommand{\Rep}{{\operatorname{Rep}}}
\newcommand{\Sym}{\operatorname{Sym}}
\newcommand{\SSym}{\operatorname{SuperSym}}
\newcommand{\Z}{\mathbb Z}
\newcommand{\ldb}{\mathopen{\{\!\!\{}} 
\newcommand{\rdb}{\mathclose{\}\!\!\}}}

\newcommand{\m}{\mathsf m}

\begin{document}
\title{Poisson algebras and Yang-Baxter equations} \author{Travis
  Schedler}
\dedicatory{Dedicated to Kazem Mahdavi, my kind mentor
  and friend.}
\maketitle

\begin{abstract}
  We connect generalizations of Poisson algebras with the classical
  and associative Yang-Baxter equations.  In particular, we prove that
  solutions of the classical Yang-Baxter equation on a vector space
  $V$ are equivalent to ``twisted'' Poisson algebra structures on the
  tensor algebra $TV$.  Here, ``twisted'' refers to working in the
  category of graded vector spaces equipped with $S_n$ actions in
  degree $n$.  We show that the associative Yang-Baxter equation is
  similarly related to the double Poisson algebras of Van den Bergh.
  We generalize to $L_\infty$-algebras and define ``infinity''
  versions of Yang-Baxter equations and double Poisson algebras.  The
  proofs are based on the observation that $Lie$ is essentially unique
  among quadratic operads having a certain distributivity property
  over the commutative operad; we also give
  an $L_\infty$ generalization. In the appendix, we prove a
  generalized version of Schur-Weyl duality, which is related to the
  use of nonstandard $S_n$-module structures on $V^{\ot n}$.
\end{abstract}
\section{Twisted Poisson algebras and the CYBE} \label{twpacybe} 
Throughout, we will work over a
characteristic-zero field $\kk$.  The tensor algebra $TV = T_\kk V$
satisfies the following twisted-commutativity property: each graded
component $V^{\ot m}$ is equipped with an $S_m$-module structure by
permutation of components, and given homogeneous elements
$v, w \in TV$ of degrees $|v|, |w|$, we have
\begin{equation}
w \ot v = (21)^{|v|,|w|} (v \ot w),
\end{equation}
where $(21)^{|v|,|w|} \in S_{|v|+|w|}$ is the permutation of the two
blocks $\{1,\ldots,|v|\}, \{|v|+1, \ldots, |v|+|w|\}$. We thus say
that $TV$ is a \textit{twisted commutative} algebra.\footnote{The
  notion of twisted algebras is an old notion from topology dating to
  at least the 1950's; see, e.g., \cite{Bar,J,Fr,LP}. They are related to
  superalgebras and color Lie algebras \cite{RWsuper,Sctgenlie}.} 

Similarly, we may define twisted Lie algebras.  Again let
$A = \bigoplus_{m \geq 0} A_m$ together with an $S_m$ action on $A_m$
for all $m$.  A twisted Lie algebra is $A$ together with a graded
bracket $\{\,,\}: A \ot A \rightarrow A$ satisfying
\begin{gather} \label{twskew}
\{w,v\} = (21)^{|v|,|w|} \{v,w\}, \\ \label{twjac}
\{u,\{v,w\}\} + (231)^{|v|,|w|,|u|} \{v, \{w,u\}\} + (312)^{|w|,|u|,|v|} 
\{w, \{u,v\}\} = 0,
\end{gather}
where $\sigma = (i_1 i_2 \ldots i_n) \in S_n$ denotes the element
$\sigma(j) = i_j$, and given $\tau \in S_3$,
$\tau^{a,b,c} \in S_{a+b+c}$ denotes the permutation acting by
permuting the blocks
$\{1,\ldots, a\}, \{a+1, \ldots, a+b\}, \{a+b+1, \ldots, a+b+c\}$. (We
will \textbf{not} use cycle notation in this paper.)

The motivating observation of this paper is as follows: If $A = TV$
is endowed with a twisted Lie algebra structure satisfying the Leibniz
rule,
\begin{equation} \label{twleib}
\{u \ot v, w\} = u \ot \{v,w\} + (213)^{|v|,|u|,|w|}  (v \ot \{u,w\}),
\end{equation}
then the Jacobi identity restricted to degree one,
$V \ot V \ot V \rightarrow T^3 V$, says that the bracket yields a skew
solution of the well-known \textbf{classical Yang-Baxter equation
  (CYBE)}:\footnote{The CYBE is a central equation in physics and the
  study of quantum groups.}
interpreted as a map $r: V \ot V \rightarrow V \ot V$, we have
\begin{gather} 
r = -r^{21}, \label{rskew} \\ \label{cybe}
[r^{12}, r^{13}] - [r^{23}, r^{12}] + [r^{13}, r^{23}] = 0,
\end{gather}
where here $r^{ij}$ denotes $r$ acting in the $i$-th and $j$-th
components (e.g., $r^{23} = \Id_V \ot r$).  The starting point for
this paper is then
\begin{thm} \label{thm1}
Let $V$ be any vector space.  Skew solutions $r \in \End(V \ot V)$ 
of the CYBE are equivalent to twisted Poisson algebra structures
on $TV$, equipped with its usual twisted commutative multiplication $\ot$.
\end{thm}
Here, a \textit{twisted Poisson} structure on $TV$ is the same as a
twisted Lie algebra structure satisfying \eqref{twleib}.

The proof is based on the twisted generalization of the following
well-known fact: a Poisson algebra structure on $\Sym V$ is the same
as a Lie algebra structure on $V$ (Proposition \ref{liepoisspr}).
Precisely, recall that an $\SSS$-module is a graded vector space
$V = \bigoplus_{m \geq 0} V_m$ together with $S_m$-actions on each
$V_m$.  $\SSS$-modules form a symmetric monoidal category, and the
notion of $\Sym V$ (the free commutative monoid in the category of
$\SSS$-modules) makes sense, and yields a twisted commutative algebra.
In the case that $V$ is concentrated in degree zero, the twisted
commutative algebra $\Sym V$, viewed as an ordinary vector space with
the induced multiplication map,\footnote{More conceptually, the forgetful
  functor from $\SSS$-modules to vector spaces is a monoidal, although
  not symmetric monoidal, functor, which is why twisted (commutative
  or associative) algebras may also be viewed as ordinary associative
  algebras.  These observations have been carried much further in,
  e.g., \cite{St,PR}.}
is the usual symmetric algebra $\Sym V_0$.  In the case $V$ is
concentrated in degree one, $\Sym V$, viewed as an ordinary vector
space with an associative multiplication, is the usual tensor algebra
$T V_1$.

Then, as explained in \S \ref{pfthm1s} below, a standard proof that
Poisson structures on $\Sym V$ are the same as Lie algebra structures
on $V$ carries over to the twisted setting, and yields Theorem
\ref{thm1}.

\begin{rem} P. Etingof pointed out to the author a connection with the
  Lie algebra $\mathfrak{tr}$ from \cite{BEER} (which is generated by
  $r_{ij}$ subject to the universal relations satisfied by $r^{ij}$
  for any skew solution $r$ of the CYBE).  More precisely, the universal
  enveloping algebra of $\mathfrak{tr}$ contains the space of all
  possible operations $V^{\ot m} \rightarrow V^{\ot m}$ obtainable from the
  twisted Lie structure on $TV$ (in terms of an indeterminate $r$).
\end{rem}

\subsection{Proof of Theorem \ref{thm1}}\label{pfthm1s}
We recall first the definition of the symmetric monoidal structure on
the category of $\SSS$-modules: Given $\SSS$-modules
$V = \bigoplus_{n \geq 0} V_n$ and $W = \bigoplus_{n \geq 0} W_n$,
\begin{equation}
(V \ot W)_{p} := 
\bigoplus_{m+n=p} \Ind_{S_m \times S_n}^{S_p} (V_m \ot W_n).
\end{equation}
As before, let $\Sym V$ denote the free commutative monoid generated
by $V$, in this category.  Theorem \ref{thm1} will follow from the
following more general result:
\begin{prop} \label{liepoisspr}
  A multiplication $\{\,,\}: V \ot V \rightarrow \Sym V$ satisfying
  \eqref{twskew}, \eqref{twjac} 
  extends uniquely to a twisted Poisson structure on $\Sym V$.
\end{prop}
In the proof below, it is helpful to have in mind the usual case when
$V$ is concentrated in degree zero; here the proof is one of the most
obvious ones of the well-known fact that a Lie structure on $V$ extends
uniquely to a Poisson structure on $\Sym V$.

We will need to introduce the following notation for technical
convenience (and it is not needed in the case $V$ is concentrated in degree
zero):
\begin{ntn}\label{sigpntn}
  Given any product of operations applied to symbols
  $x_1, \ldots, x_m$, which represent elements of degrees
  $|x_1|, \ldots, |x_m|$, let
  $\sigma'_{x_1,\ldots,x_m} \in S_{|x_1| + \ldots + |x_m|}$ be the
  permutation which corresponds to rearranging the symbols in the
  order $x_1, \ldots, x_m$ via a permutation of blocks of sizes
  $|x_1|, \ldots, |x_m|$.  For example, in the case
  $|x_1| = 2, |x_2| = 3$, we have
  $\sigma' \{x_3, x_2\} = (34512) \{x_3, x_2\}$. Also, let us allow
  $\sigma'$ to be extended linearly to linear combinations of such
  expressions.
\end{ntn}
\begin{proof} Uniqueness follows inductively from \eqref{twleib}. One
obtains the formula
\begin{equation} \label{brfla}
  \{v_1 \cdots v_m, w_1 \cdots w_n\} = \sum_{i,j} \sigma'_{v_1,\ldots,v_m,w_1,\ldots,w_n} (\{v_i, w_j\} v_1 v_2 \cdots \hat v_{i} \cdots v_m w_1 w_2 \cdots \hat w_j \cdots w_n).
\end{equation}
For existence, it suffices to verify the skew-symmetry and Jacobi identity
conditions for \eqref{brfla}.  Skew-symmetry is obvious, so it remains
to verify the Jacobi identity. This follows inductively from the
following computation:
\begin{multline}
\sigma'_{a,b,c,d} \bigl( \{ab, \{c, d\}\} + \{c, \{d, ab\}\} + \{d, \{ab, c\}\}
\bigr) \\ = \sigma'_{a,b,c,d} \bigl( b ( \{a, \{c,d\}\} + \{c, \{d, a\}\} + \{d,\{a,c\}\}) + a (\{b, \{c,d\}\} + \{c,\{d,b\}\} + \{d,\{b,c\}\}) \bigr). \qedhere
\end{multline}
\end{proof}
\begin{proof}[Proof of Theorem \ref{thm1}] Let $W$ be a vector space
  and $V := W[1]$ the associated $\SSS$-module concentrated in degree
  one. We showed in Section \ref{twpacybe} that the condition that the
  map $W \ot W \rightarrow \Sym V \cong TW$ satisfy the twisted
  skew-symmetry and Jacobi identities is exactly the statement that
  the associated element $r \in \End(W \ot W)$ is a skew solution of
  the CYBE. By Proposition \ref{liepoisspr}, we see that under this
  condition, there is a unique extension to a twisted Poisson
  structure on $\Sym V \cong TW$. Conversely, any twisted Poisson
  structure on $TW \cong \Sym V$ restricts to a skew solution
  $r \in \End(W \ot W)$ of the CYBE.
\end{proof}

\section{Double Poisson algebras} \label{dpass} The author first came
upon the aforementioned observations after reading Van den Bergh's
paper \cite{VdB} on double Poisson algebras.  These algebras formalize
Poisson geometry for noncommutative algebras such as path algebras of
quivers (see Example \ref{pathalgex} below). They are defined by the
following axioms, which are quite similar to those for twisted Lie
algebras:
\begin{defn} \cite{VdB}
A double Poisson algebra is an associative algebra $A$ with a $\kk$-linear map
$\ldb \, \rdb: A \ot A \rightarrow A \ot A$ satisfying:
\begin{gather} \label{dbskew}
\ldb a, b \rdb = - (21) \ldb b, a \rdb, \\
\label{dbjac} \sum_{i=0}^{2} (231)^i \circ \ldb - , \ldb -,
- \rdb \rdb \circ (231)^{-i} = 0, \\
\label{dbpoiss} \ldb a, bc \rdb = (b \ot 1) \ldb a, c \rdb + \ldb a, b \rdb
(1 \ot c).
\end{gather}
\end{defn}
Dropping the Poisson condition, we define:
\begin{defn}
Let $V$ be any $\kk$-vector space.
A double Lie bracket is a $\kk$-linear map
$\ldb \, \rdb: V \ot V \rightarrow V \ot V$ satisfying \eqref{dbskew}
and \eqref{dbjac}.
\end{defn}
We prove that double Lie algebras are the same as solutions of the
associative Yang-Baxter equation (AYBE), which was introduced in \cite{Agu1},\cite{Agu2} and independently in \cite{Polaybe}:
\begin{equation} \label{aybe}
r^{12} r^{13} - r^{23} r^{12} + r^{13} r^{23} = 0.
\end{equation}
Note that, when $r$ is skew ($r = -r^{21}$), then the AYBE implies the CYBE (this is an special case of \cite{Agu2}[Theorem 3.5]). Namely, let $CYBE(r)$ denote
the LHS of \eqref{cybe} and let $AYBE(r)$ denote the LHS of \eqref{aybe}.
Then, if $r$ is skew,
we have
\begin{equation} \label{cae}
CYBE(r) = AYBE(r) - (132) \circ AYBE(r) \circ (132).
\end{equation}
\begin{thm}\label{thmdpa}
  \begin{enumerate}
  \item[(i)] Let $V$ be a vector space. Double Lie algebra structures
    on $V$ are equivalent to skew solutions $r \in \End(V \ot V)$ of
    the AYBE. Hence, any double Lie algebra $V$ yields a twisted
    Poisson algebra structure on $T_\kk V$ using its tensor
    product multiplication;
\item[(ii)]  Let $A$ be an associative algebra. Suppose that
  $r \in \End(A \ot A)$ satisfies \eqref{dbpoiss}.
  If, furthermore, $r$ satisfies the CYBE, then, letting
$AYBE(r)$ denote the LHS of \eqref{aybe}, one has
\begin{equation} \label{almcybe}
(a \ot 1 \ot 1) AYBE(r) = (1 \ot 1 \ot a) AYBE(r), \quad \forall a \in A.
\end{equation}
\item[(iii)] In particular, if $A$ is a prime and noncommutative
  associative algebra, then \eqref{almcybe} implies that $AYBE(r)=0$,
  so twisted Poisson structures on $T_\kk A$ satisfying
  \eqref{dbpoiss} (where multiplication is taken in $A$) are
  equivalent to double Poisson structures on $A$.
\end{enumerate}
\end{thm}
In part (iii), ``prime'' means that, for all nonzero $a, b \in A$,
there exists $c \in A$ such that $acb \neq 0$. This is a standard
generalization of integral domains to noncommutative rings.

\begin{exam}\label{pathalgex} Consider any quiver (= directed graph) $Q$.  We recall that the \emph{path algebra}, $\kk Q$, of $Q$, is
  the algebra which, as a $\kk$-vector space, is the set of
  $\kk$-linear combinations of paths in the graph, and whose
  multiplication is given by concatenation of paths.  To be explicit,
  we may say that, for paths $p$ and $q$, $pq$ is the concatenation if
  the terminal vertex of $p$ equals the initial vertex of $q$, and
  otherwise $pq = 0$. This multiplication is extended $\kk$-linearly to
  all of $\kk Q$.  If $Q$ is \emph{strongly connected}, which means that
  for any two vertices $i$ and $j$, there is a path from $i$ to $j$, then
  $\kk Q$ is prime.  If $Q$ additionally has either at least two vertices or
  or at least two edges, then $\kk Q$ is noncommutative.  In this case, the
  theorem shows that double Poisson structures on $\kk Q$ are equivalent
  to twisted Poisson structures on $T_\kk \kk Q$ satisfying \eqref{dbpoiss}.
\end{exam}

We will say that a quiver $Q$ is (extended) Dynkin if the underlying
undirected graph (forgetting orientations, but remembering
multiplicities) is (extended) Dynkin of type $A, D$, or $E$.  
\begin{exam} For any quiver $Q$, another important algebra is
called the \emph{preprojective algebra} of $Q$, whose definition we recall as
follows.
Let $\dq \supset Q$ be the double quiver, which is the quiver with the
same set of vertices as $Q$, but with twice as many edges:
for each edge $e \in Q$, we include not merely $e$, but add an edge 
$e^* \in \dq$ which has the
same endpoints as $e$ but points in the opposite direction.  Then,
$\Pi_Q$ is defined by $\Pi_Q := \kk Q / (\sum_{e \in Q} e e^* - e^* e)$. 

Then, by  \cite[Proposition 9.2.23]{S}, for any non-Dynkin quiver $Q$,
 $\Pi_Q$ is prime, and provided $Q \neq \tilde A_0$,
  it is clear that $\Pi_Q$ is noncommutative. Thus, the theorem applies
also to this case.
\end{exam}
\begin{exam}
The deformed preprojective algebra $\Pi^\lambda_Q := \kk Q / (\lambda - \sum_{e \in Q} e e^* - e^* e)$ is filtered by degree, and its associated graded
is $\Pi_Q$.  Hence, it is also prime and noncommutative when $\Pi_Q$ is, and
thus the theorem applies.
\end{exam}
\begin{proof}[Proof of Theorem \ref{thmdpa}]
(i) Using the obvious correspondence between elements $r \in \End(V \ot V)$
and double brackets $V \ot V \rightarrow V \ot V$, the skew-symmetry condition \eqref{dbskew} becomes the condition that $r$ is skew.
Then, \eqref{dbjac} becomes
\begin{equation}
r^{12} r^{23} + r^{23} r^{31} + r^{31} r^{12} = 0.
\end{equation}
If we permute the first and third components, multiply by $-1$, and
apply skew-symmetry, we get \eqref{aybe}. This proves the first
statement. The second statement then follows from the aforementioned
fact that the AYBE implies the CYBE for skew elements $r$.

(ii) Let us write
\begin{equation}
AYBE'(r) = r^{13} r^{12} - r^{12} r^{23} + r^{23} r^{13},
\end{equation}
so that $CYBE(r) = AYBE(r) - AYBE'(r)$.  Then, using the derivation property
for $r$, one may verify that
\begin{multline}
CYBE(r) (a \ot (b_1 b_2) \ot c) = (b_1 \ot 1 \ot 1) AYBE(r) (a \ot b_2 \ot c) \\ - (1 \ot 1 \ot b_1) AYBE'(r) (a \ot b_2 \ot c) + CYBE(r)(a \ot b_1 \ot c) \cdot (1 \ot b_2 \ot 1),
\end{multline}
so that, if $CYBE(r) = 0$, then
\begin{equation}
(a \ot 1 \ot 1) AYBE(r) = (1 \ot 1 \ot a) AYBE'(r), \quad \forall a \in A.
\end{equation}
However, since $CYBE(r) = 0$, one also has $AYBE(r)=AYBE'(r)$, so
\begin{equation} \label{aybeas}
(a \ot 1 \ot 1) AYBE(r) = (1 \ot 1 \ot a) AYBE(r), \quad \forall a \in A.
\end{equation}

(iii) Assume that \eqref{aybeas} holds. Then,
\begin{equation} \label{aybe3}
(ab \ot 1 \ot 1) AYBE(r) = (a \ot 1 \ot b) AYBE(r) = (1 \ot 1 \ot ba) AYBE(r) =
(ba \ot 1 \ot 1) AYBE(r).
\end{equation}
We deduce that $((ab-ba) \ot 1 \ot 1) AYBE(r) = 0$.  So, if
$[A,A] x = 0$ implies $x =0$ for all $x \in A$, then $AYBE(r) = 0$.
This follows because, for arbitrary $y_1, y_2, y_3 \in A$,
if we write $AYBE(r)(y_1 \ot y_2 \ot y_3) = \sum_{i} v_i \ot v_i'$,
where the $v_i' \in A \ot A$ are all linearly independent, and $v_i \in A$,
then $((ab - ba) \ot 1 \ot 1) AYBE(r) = 0$ implies that
$(ab - ba) v_i = 0$ for all $i$ and all $a,b$. Hence, if $[A,A] x = 0$
implies $x=0$ for all $x$, then $v_i = 0$ for all $i$, and hence $AYBE(r)(y_1 \ot y_2 \ot y_3) = 0$. Since $y_1, y_2, y_3 \in A$ were arbitrary, $AYBE(r) = 0$
as well.

On the other hand, to say that $[A,A] a=0$ implies $a = 0$ is the same
as saying that the left ideal generated by $[A,A]$ annihilates only
zero.  But, the left ideal generated by $[A,A]$ is a two-sided ideal:
$w[x,y]z = zw[x,y] + [w[x,y],z]$.  So, for $[A,A] a = 0$ to imply
$a = 0$, it is enough that $([A,A]) a = 0$ implies $a = 0$. If $A$ is
prime and noncommutative, then $[A,A] A a = 0$ implies $a = 0$, hence
$([A,A]) a = 0$ implies $a=0$, and so $AYBE(r) = 0$, as desired.
\end{proof}

\begin{rem}
  One may be curious what happens if $A$ is commutative (such as
  $\Pi_{\tilde A_0}$).  In this case, any 
  $r \in \End(A \ot A)$ satisfying \eqref{dbpoiss}
  satisfies
\begin{equation} \label{comder1}
(a \ot 1 - 1 \ot a) r(b \ot c) = (c \ot 1 - 1 \ot c) r(b \ot a) = (b \ot 1 - 1 \ot b) r(c \ot a) = (a \ot 1 - 1 \ot a)r(c \ot b),
\end{equation}
and moreover, for any $a_1, a_2 \in A$,
\begin{equation} \label{comder2}
(a_1 \ot 1 - 1 \ot a_1)(a_2 \ot 1 - 1 \ot a_2) r(b \ot c) = (b \ot 1 - 1 \ot b)
(c \ot 1 - 1 \ot c) r(a_1 \ot a_2).
\end{equation}
So this puts special restrictions on $r$.  For instance, if $A$ is a
polynomial algebra over a field, then so is $A \ot A$, and by unique
factorization, \eqref{comder2} implies that
$r(b \ot c) = (b \ot 1 - 1 \ot b)(c \ot 1 - 1 \ot c) f = r(c \ot b)$ for all
$b,c \in A$ and some fixed $f$ in the quotient field of $A \ot A$.
(Furthermore, one must have $f \in A \ot A$ unless $A$ is a polynomial
algebra in only one variable $x$, in which case
$f \in (A \ot A)\cdot (x \ot 1 - 1 \ot x)^{-2}$.)  Also, for $r$ to be
skew, $f$ must be skew, and in this case, the image of $r$ lies in
$A \wedge A$, so one deduces that
$r^{12} r^{13} = -(213) \circ (r^{12} r^{13}) = r^{12} r^{23}
\circ (213)$.
In the one-variable case, one deduces that $CYBE(r) (x \ot x \ot x) = 0$
iff $AYBE(r) (x \ot x \ot x) = 0$, and the latter is true (using the
Poisson condition) iff $AYBE(r) = 0$ (and hence $CYBE(r)=0$).  One may
then deduce inductively that these are all zero iff
$r(x \ot x) = \lambda (x \ot 1 - 1 \ot x)$ for some $\lambda \in \kk$.  It
would be interesting to see if there are other solutions in
more variables (and whether AYBE and CYBE have the same solutions).
Note that $\Pi_{\tilde A_0}$ is the two-variable case.
\end{rem}
\begin{rem} Many solutions of the AYBE with (graded and spectral)
  parameters $u,v$, related to solutions of the CYBE and QYBE, have
  been classified in \cite{Polaybe,Strig,Polaybe2}.  For example, any
  ``associative Belavin-Drinfeld structure'' gives rise to a
  (trigonometric) such solution which additionally satisfies the QYBE
  and CYBE with spectral parameters.  One may interpret a solution of
  the AYBE and skew-symmetry (unitarity) with parameters as a sort of
  graded version of double Lie algebra; to make it Poisson, one would
  need to find a compatible multiplication (if it exists).
\end{rem}

\section{Operadic generalization} \label{opsec} We ask the
question: Is Theorem \ref{thm1} a special property of Lie algebras, or
can it be generalized to other operads (suitably replacing the CYBE
with other equations)? To make sense of this question, for any 
twisted-commutative algebra $A$ equipped with an additional binary
operation $\star: A
\ot A \rightarrow A$, we generalize the Leibniz rule \eqref{twleib} to
\begin{equation}\label{genleib}
(uv) \star w = u (v \star w) + (213)^{|v|,|u|,|w|} v (u \star w).
\end{equation}
Let $\cO$ be any operad generated by a single element $\m \in \cO(2)$.
If $A$ as above is an $\cO$-algebra (with $a \star b := \m(a \ot b)$)
satisfying \eqref{genleib}, we call it a \emph{twisted distributive
  $\cO$-algebra}.

Let $\F$ be the operad freely generated by 
an element also denoted by $\m \in \F(2)$.
By a quadratic operad $\cO = \F / (R)$, we mean one such that $R = R_2
\oplus R_3$, with $R_2 \subset \F(2)$ and $R_3 \subset \F(3)$.  In
particular, the only possible relations from $R_2$ are symmetry
conditions on $\m$, namely, that $(21) \m = \pm \m$ as elements of
$\cO$. Also, $R_3$ consists of relations which are quadratic in $\m$.
We then prove the following result
(in \S \ref{thm2s}):
\begin{defn} Let $\text{\textit{Lie-adm}}$ denote the Lie-admissable operad, which is
  the operad whose ordinary algebras are vector spaces together with a
  binary operation whose skew-symmetrization is a Lie bracket.
\end{defn}
\begin{thm} \label{thm2}
  The only
  quadratic operads $\cO = \F / (R)$ for which distributive
  $\cO$-algebra structures on $\Sym V$ are equivalent to $\cO$-algebra
  structures on $V$,
 for  every vector space $V$, are the following five:
\begin{itemize}
\item $\cO = Lie$, the Lie operad
\item $\cO = \text{\textit{Lie-adm}}$, the Lie-admissable operad
\item $\cO = \F / (\m + (21) \m)$, the operad whose algebras are
vector spaces with a skew-symmetric binary operation
\item $\cO = \F / (\m - (21) \m)$, the operad whose algebras are
vector spaces with a symmetric binary operation
\item $\cO = \F$, the free operad generated by $\m$, the operad
whose algebras are magmas.
\end{itemize}
These are also exactly the operads for which twisted distributive
$\cO$-algebra structures on $\Sym V$ are equivalent to operations
$V \ot V \rightarrow \Sym V$ satisfying the relations $R$, for every
$\SSS$-module $V$.
\end{thm}
After proving this theorem, we will generalize to quadratic operads
which are not generated by only a single operation, and find that
$L_\infty$-algebras (heuristically, Lie algebras up to homotopy) are
the prototypical example of algebras with operations satisfying the
desired distributivity property (Theorem \ref{thm3} and Remark
\ref{thm3rem}).
\subsection{Proof of Theorem
  \ref{thm2}} \label{thm2s} As in Proposition \ref{liepoisspr}, it is
not difficult to show that the listed operads have the desired
property.  We show the converse.  Let $\F$ be as in Theorem
\ref{thm2s}, and let $\cO = \F / (R)$ with $R = R_2 \oplus R_3$, where
$R_2 \subset \F(2)$ and $R_3 \subset \F(3)$.  We show that, if
$\cO$-algebra structures on every vector spaces $V$ are equivalent to
distributive $\cO$-algebra structures on $\Sym V$ , then $\cO$ is
one of the listed operads.

Assume first that $R_2 = 0$, i.e., there is no (skew)-symmetry axiom for
$\cO$-algebras. Take an arbitrary element of $R_3$: this is equivalent
to a quadratic axiom for $\cO$-algebras $(V, \star)$. Let us write it as
\begin{equation} \label{arbreln} \sum_{\sigma \in S_3}
  \lambda_{\sigma, 1} b_{\sigma(1)} \star (b_{\sigma(2)} \star
  b_{\sigma(3)}) + \lambda_{\sigma, 2} (b_{\sigma(1)} \star
  b_{\sigma(2)}) \star b_{\sigma(3)} = 0,
\end{equation}
for some constants $\lambda_{\sigma, i}$.

If $\cO$-algebra structures $(V, \star)$ on $V$ are equivalent to
distributive $\cO$-algebra structures on $\Sym V$ for all $V$, then if
we expand the above axiom for $(b_1, b_2, b_3) = (b_1' b_1'', b_2,
b_3)$ using the Leibniz rule \eqref{genleib}, the terms of the form
$(x \star y) (z \star w)$ must cancel identically, for $\{x,y,z,w\} =
\{b_1',b_1'',b_2,b_3\}$.  That is, the following expression must be identically zero:
\begin{multline}
\sum_{\sigma \mid \sigma(1) = 1} \lambda_{\sigma,2} ((b_1' \star b_2) (b_1'' \star b_3) + (b_1' \star b_3) (b_1'' \star b_2)) \\ + \sum_{\sigma \mid \sigma(2) = 1} (\lambda_{\sigma, 1} + \lambda_{\sigma,2}) ((b_{\sigma(2)} \star b_1') (b_1'' \star b_{\sigma(3)}) + (b_{\sigma(2)} \star b_1'') (b_1' \star b_{\sigma(3)}) ) + \sum_{\sigma \mid \sigma(3) = 1} \lambda_{\sigma, 1} ((b_2 \star b_1') (b_3 \star b_1'') + (b_3 \star b_1')(b_2 \star b_1'')).
\end{multline}
This can only happen if the following equations are satisfied:
\begin{equation} \label{reqreln1}
\lambda_{\id, 2} = - \lambda_{(132), 2}, \quad \lambda_{(213), 2} = - \lambda_{(213), 1}, \quad \lambda_{(312), 2} = - \lambda_{(312), 1}, \quad \lambda_{(231), 1} = - \lambda_{(321),1}.
\end{equation}
Similarly, if we plug in instead $(b_1, b_2, b_3) = (b_1, b_2' b_2'', b_3)$ or $(b_1, b_2, b_3) = (b_1, b_2, b_3' b_3'')$, we obtain additionally the
following conditions:
\begin{gather} 
\lambda_{(213), 2} = -\lambda_{(231), 2}, \quad \lambda_{\id, 2} = - \lambda_{(321),1}, \quad \lambda_{(321), 2} = - \lambda_{\id,2}, \quad \lambda_{(132), 1} = - \lambda_{(312), 1}, \\
\lambda_{(312), 2} = - \lambda_{(321), 2}, \quad \lambda_{(132), 2} = - \lambda_{(231), 1}, \quad \lambda_{(231),2} = - \lambda_{(132),1}, \quad \lambda_{\id, 1} = - \lambda_{(213), 1}.
\end{gather}
(These can also be obtained from \eqref{reqreln1} by applying the action of $S_3$.)

We deduce that the only possible element of $R_3$ is a multiple of the associated Lie relation, which proves the theorem in the case that $R_2 = 0$.

If $R_2 \neq 0$, then the above computation simplifies. Suppose that
$R$ is spanned by $\m - \varepsilon (21) \m$ for $\varepsilon \in \{1, -1\}$.
Take an arbitrary element of $R_3$ and write the corresponding axiom for
$\cO$-algebras as
\begin{equation}\label{arbrelnsymm}
 \lambda_1 b_1 \star (b_2 \star b_3) + \lambda_2 b_2 \star (b_3 \star b_1) + \lambda_3 b_3 \star (b_1 \star b_2) = 0.
\end{equation}
Next, we plug in $b_1 = b_1' b_1''$ and gather all terms on the LHS of the form
$(x \star y)(z \star w)$ where $\{x,y,z,w\} = \{b_1',b_1'', b_2, b_3\}$:
\begin{equation}\label{leibcalcsymm}
\lambda_2 \bigl( (b_2 \star b_1') (b_3 \star b_1'') + (b_2 \star b_1'') (b_3 \star b_1') \bigr) + \lambda_3 \bigl( (b_3 \star b_1') (b_1'' \star b_2) + (b_3 \star b_1'') (b_1' \star b_2) \bigr).
\end{equation}
The above must be zero, in order for $\cO$ to have the desired property.
Using the symmetry condition $x \star y = \varepsilon y \star x$, we may rewrite \eqref{leibcalcsymm} as
\begin{equation}\label{leibcalc2symm}
(\lambda_2 + \varepsilon \lambda_3) \bigl( (b_2 \star b_1') (b_3 \star b_1'')
+ (b_3 \star b_1') (b_2 \star b_1'') \bigr).
\end{equation}
For this to be zero in general, we require that $\lambda_2 = -
\varepsilon \lambda_3$.  Similarly, setting $(b_1, b_2, b_3) = (b_1,
b_2' b_2'', b_3)$, we conclude that $\lambda_1 = - \varepsilon
\lambda_3$, and finally we conclude that $\lambda_1 = - \varepsilon
\lambda_2$.  Hence, if $\varepsilon = 1$, then $\lambda_1 = \lambda_2
= \lambda_3 = 0$, and if $\varepsilon = -1$, then $\lambda_1 =
\lambda_2 = \lambda_3$ can be arbitrary.  This proves the desired result.

For the final statement, it suffices to generalize Proposition
\ref{liepoisspr} to the listed operads. This is straightforward and
omitted. \qed

\subsection{$L_\infty$ generalization and arbitrary quadratic operads}
In this subsection we will drop the assumption on $\cO$ that it be
generated by a single binary operation.  Suppose instead that $\cO$ is
generated by any operations, of any arity, satisfying quadratic
and linear relations.  

The simplest example of this is the case where $\cO$
includes a differential $d \in \cO(1)$, satisfying $d^2 = 0$.  We see
that, already, it is not true that $d^2|_{V} = 0$ implies $d^2|_{\Sym V} = 0$.  The solution to this problem is well-known: make $V$ a graded vector space,
take $d$ to be an operator of degree $-1$, and make $\Sym V$ the supersymmetric algebra generated by $V$, i.e., $v w = (-1)^{|v| |w|} wv$.

It then turns out that Theorems \ref{thm1}, \ref{thm2} generalize,
roughly, by replacing Lie algebras with $L_\infty$ algebras
(heuristically, these are Lie algebra up to homotopy).

We note that the ``super'' grading above is independent of the twisted
grading, so that when one has both, $V$ is bigraded.  To simplify
things, we restrict to ordinary (not twisted) algebras with a single
grading. All of the results generalize to the twisted setting, by
working with $\SSS$-modules and adding permutations $\sigma'$ (using
Notation \ref{sigpntn}) to the beginning of formulas.

Recall that an $L_\infty$ algebra is a $\Z$-graded vector space $A$
equipped with a differential $d = \{\ \}_1: A \rightarrow A$ of degree
$1$, and completely graded-skew-symmetric operations
$\{\ \}_n: A^{\ot n} \rightarrow A$ of degrees $2-n$ (for
$n \geq 2$), satisfying axioms
\begin{equation} \label{linfax}
\sum_{i+j = m+1, i,j \geq 1} \sum_{\sigma \in S_{m}} (-1)^{i} \sgn_{odd}(\mathbf{a},\sigma) \{\{a_{\sigma(1)}, \ldots, a_{\sigma(i)}\}_i,a_{\sigma(i+1)},\ldots, a_{\sigma(i+j-1)}\}_j = 0,
\end{equation}
for all $m\geq 1$, where $\sgn_{odd}(\mathbf{a},\sigma)$ is
the sign of the permutation $\sigma^{|a_1|+1,|a_2|+1, \ldots, |a_{m}|+1} \in
S_{|a_1|+1 + |a_2|+1 + \cdots + |a_{m}|+1}$ obtained from $\sigma$ by
acting on blocks of the sizes $|a_1|+1, \ldots, |a_{m}|+1$.
  
In particular, this includes the axiom that $d$ is a differential, and
that the Jacobi identity for $\{\ \}_2$ is satisfied up to chain homotopy.

We will use the notation
\begin{equation}\label{bntn}
\mathbf{b} := (a_1, a_2, \ldots, a_{k-1}, a_k', a_k'', a_{k+1}, \ldots, a_m).
\end{equation}
We define the Leibniz rule for an $L_\infty$ algebra endowed with an additional \textbf{super}commutative multiplication as
\begin{multline} \label{gennerleib}
\{a_1, a_2, \ldots, a_{k-1}, a_k' a_k'', a_{k+1}, \ldots, a_m\}_m \\ = 
(-1)^{|a_k'| (|a_1| + |a_2| + \cdots + |a_{k-1}|)} 
a_k' \{a_1, a_2, \ldots, a_{k-1}, a_k'', a_{k+1}, \ldots, a_m\}_m
\\ + (-1)^{|a_k''| (|a_1| + |a_2| + \cdots + |a_{k-1}|+|a_k'|)} a_k'' \{a_1, a_2, \ldots, a_{k-1}, a_k', a_{k+1}, \ldots, a_m\}_m.
\end{multline} 
\begin{thm} \label{thm3}
An $L_\infty$ structure on $\SSym V$ satisfying the Leibniz rule \eqref{gennerleib}
is equivalent to operations 
$\{\ \}_i: V^{\ot i} \rightarrow V$ satisfying the $L_\infty$ axioms
\eqref{linfax}.
\end{thm}
\begin{proof}
  This is similar to the proof of Proposition \ref{liepoisspr}. We
  need to show that, given an operation $V \ot V \rightarrow \SSym V$
  satisfying the $L_\infty$ axioms \eqref{linfax}, there is a unique
  extension using the Leibniz rule \eqref{gennerleib} to a
  multiplication $\SSym V \ot \SSym V \rightarrow \SSym V$, and this
  also satisfies the $L_\infty$ axioms.  The fact that there is a
  unique extension is easy.  To show it satisfies the $L_\infty$
  axioms, it is enough inductively to verify that, for any $m$-tuple
  of the form $(a_1, a_2, \ldots, a_m) = (a_1' a_1'', a_2, a_3, \ldots, a_m) \in (\SSym V)^m$, then the
  $L_\infty$-axioms for $(a_1', a_2, a_3, \ldots, a_m)$ and $(a_1'', a_2, a_3, \ldots, a_m)$ imply the $L_\infty$ axioms for $(a_1, a_2, \ldots,
  a_m)$. This follows by expanding the expression, for all $m \geq 1$:
\begin{equation} \label{linfax2}
\sum_{i,j: i+j=m+1} \sum_{\sigma \in S_m} (-1)^i \sgn(\mathbf{a}, \sigma) \{\{a_{\sigma(1)}, \ldots, a_{\sigma(i)}\}_i,a_{\sigma(i+1)},\ldots, a_{\sigma(m-1)}\}_j 
\end{equation}
using \eqref{gennerleib}, and verifying that the terms of the form
$\pm \{ \ldots \}_i\cdot \{ \ldots \}_j$ cancel. In more detail, it is
equivalent to sum above not over all permutations $\sigma \in S_m$,
but only over the $i,j-1$-shuffles: that is, $\sigma$ such that
$\sigma(\ell) < \sigma(\ell+1)$ for all $\ell \neq i$ (these are the
permutations that leave the order of $1,2,\ldots,i$ and
$i+1,i+2,\ldots,i+j-1=m$ unchanged). Then, there is a canonical
bijection between $i,j-1$ shuffles $\sigma$ such that $\sigma(1)=1$
and $j,i-1$ shuffles $\sigma'$ such that $\sigma'(1)=1$, and the terms
of the form $\pm \{\ldots \}_i \cdot \{\ldots \}_j$ that appear in the
expansion of the summand of \eqref{linfax2} corresponding to
$(i,j,\sigma)$ cancel with the terms of the form $\pm \{\ldots \}_j
\{\ldots\}_i$ that appear in the expansion of the summand of
\eqref{linfax2} corresponding to $(j,i,\sigma')$. We omit further details.
\end{proof}
\begin{rem}\label{thm3rem}
  One may obtain a converse of the above theorem, parallel to Theorem
  \ref{thm2}, that is, a description of all $\Z/2$-graded quadratic
  operads $\cO$ with the super version of the distributivity property
  of Theorem \ref{thm2} 
  (without the
  condition that the operad be generated by a single multiplication), using
 the Leibniz rule \eqref{gennerleib} with $\{ \}_m$ replaced by
  any $m$-ary operation corresponding to a generator of $\cO$.

  Namely, suppose that the $\cO$ is a $\Z/2$-graded quadratic operad
  generated by totally graded-skew-symmetric operations $o_1, \ldots,
  o_m$ (with no linear relations aside from that the $o_i$ are
  graded-skew-symmetric). Then, the distributivity property holds iff
  the projection of the quadratic relations to the $\SSS$-module
  generated by $o_i \circ (o_j \ot \id), o_j \circ (o_i \ot \id)$
  consists at most of the relation
\begin{multline}
\sum_{\sigma \in S_{|o_i|+|o_j|-1}} \sgn_{odd}(\mathbf{a},\sigma)\bigl( (-1)^{|o_i|} \{\{a_{\sigma(1)}, \ldots, a_{\sigma(|o_i|)}\}_{i},a_{\sigma(|o_i|+1)},\ldots, a_{\sigma(|o_i|+|o_j|-1)}\}_{j}
\\ + (-1)^{|o_j|} \{\{a_{\sigma(1)}, \ldots, a_{\sigma(|o_j|)}\}_{j},a_{\sigma(|o_j|+1)},\ldots, a_{\sigma(|o_i|+|o_j|-1)}\}_{i}\bigr) = 0,
\end{multline}
where $\{\ \}_i, \{\ \}_j$ denote applying the operations $o_i, o_j$
to the given arguments.  The proof is similar to the proofs of
Theorems \ref{thm2} and \ref{thm3}.  If $\cO$ is not generated by
totally skew-symmetric operations, then the only allowable quadratic
relations are those specifying that the skew-symmetrization of the
generating operations satisfy certain relations as above.  We omit the
details.
\end{rem}

\section{Yang-Baxter-infinity equations and double Poisson-infinity
  algebras}
In view of the fact that $L_\infty$ algebras also have the
distributivity property of Theorem \ref{thm2}, we explain here the double
Poisson analogue of twisted distributive $L_\infty$ structures on
$TV$, which we call ``double Poisson-infinity algebras'' (Definition
\ref{dpinfd}).  Here, ``infinity'' refers to relaxing the Jacobi
identity up to higher homotopies: double Poisson-infinity algebras, as
we define them, still include an honest associative algebra and the
bracket is still skew-symmetric.

Further, using twisted and double Poisson-infinity algebras, we define
infinity versions of the classical and associative Yang-Baxter
equations (Definition \ref{ybeinfd}), by analogy with Theorem
\ref{thm1}. The $\mathrm{CYBE}_\infty$ yields equations for
\textit{sequences} of elements $r_n \in \mathfrak{g}^{\ot n}$ where
$\mathfrak{g}$ is any graded Lie algebra, and the
$\mathrm{AYBE}_\infty$ yields equations for $r_n \in A^{\ot n}$, where
$A$ is any graded associative algebra. We do not know if there exists
a corresponding notion of quantum Yang-Baxter equation-infinity.

Specifically, we will (abusively) call a (twisted) commutative and
$L_\infty$ algebra satisfying \eqref{gennerleib} a (twisted)
``Poisson-infinity'' algebra.\footnote{This is abusive because this
  notion of Poisson-infinity only relaxes the Jacobi identity of the
  bracket, but not the associativity of the multiplication. The
  operadic notion of Poisson-infinity relaxes both, as well
  as the Leibniz rule.}
As an application of our comparison of twisted and double Poisson
algebras, it makes sense to define double Poisson-infinity algebras.
To do this, we need only define a ``double'' version of the
Jacobi-infinity identity \eqref{linfax}.  As in the usual setting,
we do this by replacing sums over all permutations by sums over only
cyclic permutations:
\begin{defn}\label{dpinfd}
A double Poisson-infinity algebra is a $\Z$-graded 
associative algebra $A$ together with
brackets $\{\ \}_n: A^{\ot n} \rightarrow A^{\ot n}$ of degree $2-n$, for all $n \geq 1$, satisfying the identities (for all $n \geq 1$):
\begin{gather}
\text{Skew-symmetry: } \quad \sigma \{a_{\sigma(1)}, \ldots, a_{\sigma(n)}\}_n 
= \sgn(\mathbf{a}, \sigma) \{a_1, a_2, \ldots, a_n\}, \forall \sigma \in S_n,\\
\text{Jacobi}_\infty:  \label{jinfax}
\sum_{i+j = n+1} \sum_{\sigma \in \Z/(i+j-1)} (-1)^{i} \sgn_{odd}(\mathbf{a},\sigma) \{\{a_{\sigma(1)}, \ldots, a_{\sigma(i)}\}_i,a_{\sigma(i+1)},\ldots, a_{\sigma(i+j-1)}\}_j = 0,
\end{gather}
\begin{multline}
\text{Double Leibniz: } \{a_1, a_2, \ldots, a_{n-1}, a_n' a_n''\}_n \\ = (-1)^{|a_n'| (|a_1| + \ldots + |a_{n-1}|)} a_n' \{a_1, \ldots, a_{n-1}, a_n''\}_n
+ (-1)^{n|a_n''} \{a_1, \ldots, a_{n-1}, a_n'\} a_n''. 
\end{multline}
\end{defn}
In a future paper, we hope to explain a double version of Kontsevich's
formality theorem \cite{Kform}, where the above will replace
$L_\infty$ for the Poisson side (the differential operator side
will use \cite{GS}).  

Finally, we obtain infinity versions of the Yang-Baxter equations by
writing down the Poisson-infinity conditions in terms of elements
$r_n \in \End(V^{\ot n})$ (which we may generalize to
$\mathfrak{g}^{\ot n}, A^{\ot n}$).  In order to make the sum over as
few terms as possible, and to specialize to the ordinary Yang-Baxter
equations, we use the
\begin{ntn} Let $Sh_{i,j} \subset S_{i+j}$ denote the set of $i,j$-shuffles:
this means permutations $(k_1 k_2 \cdots k_{i+j})$ such that $k_1 < k_2 < \ldots < k_i$ and $k_{i+1} < k_{i+2} < \cdots < k_{i+j}$.
\end{ntn}
\begin{defn}\label{ybeinfd} (i) Let $\mathfrak{g}$ be a graded Lie algebra.  The \textbf{classical Yang-Baxter-infinity equations} for elements $\{r_n \in \mathfrak{g}^{\ot n}\}_{n \geq 1}$  of degrees $2-n$ are, for all $n \geq 1$,
\begin{equation}
\sum_{i+j = n+1} (-1)^i
\sum_{\sigma \in Sh_{i,i+j-1}} [r_{i}^{\sigma(1),\sigma(2),\ldots,\sigma(i)}, 
r_j^{\sigma(1),\sigma(i+1),\sigma(i+2),\ldots,\sigma(i+j-1)}] = 0.
\end{equation}

(ii) Let $A$ be a graded associative algebra. The \textbf{associative Yang-Baxter-infinity equations} for elements $\{r_n \in A^{\ot n}\}_{n \geq 1}$ of degrees $2-n$ are, for all $n \geq 1$,
\begin{equation}
\sum_{i+j = n+1} (-1)^i
\sum_{\sigma \in \Z/n} r_{i}^{\sigma(1),\sigma(2),\ldots,\sigma(i)}
r_j^{\sigma(1),\sigma(i+1),\sigma(i+2),\ldots,\sigma(i+j-1)} = 0.
\end{equation}
\end{defn}
We remark that the graded condition is not essential for the
$\mathrm{AYBE}_\infty$ to make sense, although it is needed for
$\mathrm{CYBE}_\infty$ since we need supercommutators.  

It would be interesting to see if there is any reasonable
$\infty$-analogue of the quantum Yang-Baxter equation, obtained by
somehow ``quantizing'' the above equation.

\section{Non-Poisson twisted algebra structures on $TV$}
Let $V$ be a vector space. One interpretation of Theorem \ref{thm2} is
that the Leibniz rule \eqref{twleib} is not a good condition to impose
for many types of algebra structures on $TV$.  For example, the result
holds for neither associative algebras nor commutative algebras.  In
this section, in the form of remarks, we briefly explain how twisted
associative algebra structures on $TV$ (\textit{without} the Leibniz
condition) are related to the \textit{quantum} rather than the
classical Yang-Baxter equation.

Let $V$ be a representation $\rho$ of a Hopf algebra $H$, which endows
$V^{\ot m}$ (using $\rho^{\ot m}$) and hence $TV$ with a canonical
structure of $H$-representation.  We look for twisted associative
algebra structures on $TV$ given by a single element $\J \in H \ot H$,
by the rule
\begin{equation} \label{Jrep}
a \cdot b = \rho^{\ot |a|} \ot \rho^{\ot |b|}(\J) (a \ot b).
\end{equation}
We obtain the following:
\begin{prop}
The formula \eqref{Jrep} yields a twisted associative algebra structure
iff $J$ is a twist:
\begin{equation} \label{twisteq}
(\Delta \ot 1)(\J) (\J \ot 1) = (1 \ot \Delta)(\J) (1 \ot \J).
\end{equation}
\end{prop}
This has a generalization to algebras over any operad $\cO$, where $\J$
is replaced by an element satisfying an $\cO$-version of
\eqref{twisteq}.

We may generalize the above to the case where there is no universal
element $\J$.  First, note that the induced multiplication operation
on $TV$ is described by the restrictions to
$W_1 \ot W_2 \rightarrow TV$, where $W_1, W_2$ are irreducible
representations of $H$ occurring in $TV$.  This is still true without
an element $\J$.  These maps need only satisfy associativity for
triples $W_1, W_2, W_3$, and preserve a restricted form of the
permutation action.

Furthermore, we don't need to be given an $H$, since we can always
take $H = U(\mathfrak{gl}(V))$.  Let $\Rep_{alg}(\End(V))$ denote the
category of algebraic representations of $\End(V)$ viewed as a
$\kk$-algebraic monoid, which is equivalent to the category of
representations of $H$ occurring in $TV$. (This includes all
representations up to twisting by the trace representation of
$\mathfrak{gl}(V)$ (and up to isomorphism).)  By Schur-Weyl duality,
the $S_n$ action on $V^{\ot n}$ spans all of $\End_H(V^{\ot n})$.  In
the language of category theory, we obtain the
\begin{prop}
  Twisted associative algebra structures on $TV$ are the same as
  monoidal structures on the fiber functor
  $\Rep_{alg}(\End(V)) \rightarrow Vect$.
\end{prop}

The above also naturally generalizes to the case of algebra structures
on $TV$ that use a modified permutation action, given by an element
$R \in \End(V \ot V)$ (so, $R^{21} R = \Id$ and $R$ satisfies the
quantum Yang-Baxter equation,
$R^{12} R^{13} R^{23} = R^{23} R^{13} R^{12}$).  In this case, one has
a natural Hopf algebra $H_R$, defined in \cite{FRT} (see also
the appendix), which makes $V$ a
canonical comodule.  Using a generalized version of Schur-Weyl duality
(Theorem \ref{gswdt}),  we obtain the
\begin{prop}
  Twisted associative algebra structures on $(TV,R)$ are the same as
  monoidal structures on the fiber functor
  $Comod(H_R) \rightarrow Vect$.
\end{prop}

\section{Acknowledgements}
The author happened upon the initial discovery in the context of
discussions (aimed at \cite{GS}) with Victor Ginzburg, and he would
like to thank Ginzburg for those discussions as well as for
encouragement.  Thanks are also due to Pavel Etingof and Dimitri
Gurevich for useful comments, and to the anonymous referee for many
helpful suggestions which improved the exposition. Finally, the author
would like to thank K. Mahdavi, D. Koslover, the University of Texas
at Tyler, and participants of the 2007 conference in Tyler for the
opportunity to communicate this research and interesting discussions.
This research was partially supported by an NSF GRF.

\appendix
\section{Generalized Schur-Weyl duality}\label{hrsect}
\subsection{The \cite{FRT} construction}
We recall the definition of the coquasitriangular Hopf algebra $H_R$.
  Beginning with any solution $R$ of the
QYBE, $R^{12} R^{13} R^{23} = R^{23} R^{13} R^{12}$,
Faddeev, Reshetikhin, Takhtajan,
and Sklyanin \cite{FRT} constructed the following bialgebra, which is
like a $R$-twisted version of $\cO(End(V))$, the commutative bialgebra
of functions on the multiplicative monoid $\End(V)$. It has a
coquasitriangular structure from which one recovers $R$.
\begin{defn}\cite{FRT}
Define $H_R$ to be the quotient of the free algebra $F :=
\kk \langle L_{ij} \rangle_{i,j \in \{1,\ldots,n\}}$ by the
following relations.  Set
 $L = \sum_{i,j \in \{1,\ldots,n\}} e_{ij} \ot L_{ij}
 \in \End(V) \ot F$. Then $H_R := F / I_R$, where
$I_R$ is the ideal generated by the relations
\begin{equation} \label{hrrele}
R^{12} L^{13} L^{23} = L^{23} L^{13} R^{12} \quad \in \End(V) \ot \End(V) \ot F.
\end{equation}
We also let $L$ denote its own image under the quotient
$F \onto H_R$.  Then, the coproduct $\Delta$ and counit $\epsilon$
are defined by
\begin{equation}
\Delta(L) = L^{12} L^{13}, 
\quad \epsilon(L) = 1.
\end{equation}
\end{defn}
\begin{thm} \cite{FRT} (cf.~\cite{Kas})
The preceding definition makes sense and defines a bialgebra (with
a unique coquasitriangular structure inducing $R$).
\end{thm}
By ``coquasitriangular structure inducing $R$,'' we mean a map
$H_R \ot H_R \rightarrow \kk$, which satisfies the dual of the quasitriangularity axioms (replacing multiplication by convolution), and whose action on the
standard comodule is $R$. (See \cite{Kas} for details.)

We will need the standard comodule:
\begin{defn}
The ``standard comodule'' $V$ of $H_R$ is given by the element
$L \in \End(V) \ot H_R \cong \Hom(V, V\ot H_R)$.  That is, the map
$\Delta: V \rightarrow V \ot H_R$ is given by $\Delta(e_j) = \sum_{i=1}^{\text{dim }V} e_i \ot L_{ij}$. Let $V^{\ot m}$ denote the comodules obtained
by the tensor power of this one.
\end{defn}
\subsection{Generalized Schur-Weyl duality}
In this subsection we prove a generalization of Schur-Weyl duality to
the bialgebras $H_R$.  This says that the irreducible
subrepresentations of $V^{\ot m}$ are given by Young diagrams using the
$R$-symmetric action.  This result should not be too surprising, given
that the relations \eqref{hrrele} are defined in terms of $R$, and it
is possible that the result even motivated the definition of $H_R$.
However, since the author could not find it in the literature, the result
is given here.  This result shows that $Comod(H_R)$ is the
same as the category $\mathcal{SW}(V)$ studied in \cite{GurM} (under certain
conditions on $R$).
\begin{ntn} For any permutation $\sigma \in S_n$, let $\tau_\sigma: V^{\ot n} \rightarrow V^{\ot n}$ denote the permutation of components (i.e., the standard permutation action).
\end{ntn}
\begin{thm}\label{gswdt} Assume $\text{char}(\kk) = 0$ and $R$ is a unitary
  solution of the QYBE.
  Let $SR_m \subset \End(V^{\ot m})$ be the image of the $R$-symmetric
  action of $S_m$ (generated by elements $\tau_{b,b+1} R^{b,b+1}$),
  and let $HR_m \subset \End(V^{\ot m})$ be the span of linear maps of
  the form $\phi \circ \Delta$, where
  $\Delta: V^{\ot m} \rightarrow V^{\ot m} \ot H_R$ is the comodule
  action, and $\phi \in \End(H_R, \kk)$ is any linear map.  Then, one has
\begin{gather}
\End_{SR_m}(V^{\ot m}) = HR_m, \label{srmhrm} \\
\End_{HR_m}(V^{\ot m}) = SR_m, \label{hrmsrm} \\
V^{\ot m} \cong \bigoplus_{\lambda} \rho_{\lambda,S_m} \ot \rho_{\lambda,H_R},
\label{shrm}
\end{gather}
where the sum is over Young diagrams $\lambda$ parametrizing
irreducible representations of $S_m$, and $\rho_{\lambda,S_m}$ is the
corresponding irreducible representation of $S_m$. The space
$\rho_{\lambda,H_R}$ is the $H_R$-subcomodule of $V^{\ot m}$ equal to
$c_\lambda(R) V^{\ot m}$, where $c_\lambda$ (cf.~\cite{FH}) is the
Young symmetrizer in $\kk[S_m]$ corresponding to $\lambda$ (so that
$k[S_m] a_\lambda \cong \rho_{\lambda,S_m}$), and $c_\lambda(R)$ is
the corresponding element of $\End(V^{\ot m})$ given by the
``$R$-permutation action'' $*_R$ of $S_m$:
\begin{equation} \label{rpae}
(b,b+1) *_R v := \tau_{(b,b+1)} R^{b,b+1} \cdot v.
\end{equation}
Furthermore, $\rho_{\lambda,H_R}$ is an irreducible $H_R$-comodule, and
\eqref{shrm} is the multiplicity-free decomposition of $V^{\ot m}$ into
irreducible $SR_m \otimes HR_m$-modules.
\end{thm}
\begin{proof}
  We claim that \eqref{srmhrm} is true. First, note that $H_R$ is
  graded (setting $T_{i,j}$ to have degree 1), and that
  $\Delta(V^{\ot m}) \subset V^{\ot m} \ot H_R[m]$.  Next, consider
  $H_R[1]=F[1]$ to be $\End(V)^*$, by the pairing
  $(e_{i,j},T_{k,\ell}) = \delta_{i,k} \delta_{j,\ell}$.  This means
  that $H_R[m]^* \subset F[m]^* = \End(V)^{\ot m}$, which is the
  subspace respecting the relation \eqref{hrrele}.
  
  Then, the LHS of \eqref{hrrele}, $R^{12} L^{13} L^{23}$, is
  identified, as a subspace of
  $\End(V) \ot \End(V) \ot F[2] \cong \End(V) \ot \End(V) \ot
  (\End(V)^*)^{\ot 2}$,
  with $R^{12} (\Id_{\End(V)}^{13} \Id_{\End(V)}^{24})$, where
  $\Id_{\End(V)} \in \End(V) \ot \End(V)^* \cong \End(\End(V))$ is the
  canonical identity element.  Thus, for any $\phi \in \End(V)^{\ot 2}$, we have
\begin{equation}
\phi^{34} (R^{12} \Id_{\End(V)}^{13} \Id_{\End(V)}^{24})
= R \phi \in \End(V \ot V).
\end{equation}
Similarly, $\phi^{34}$ applied to the RHS of \eqref{hrrele}
(considered as an element of $\End(V)^{\ot 2} \ot (\End(V)^*)^{\ot 2}$)
is identified with $\phi^{21} R$.  So, the condition for $\phi$ to be
an element of $H_R[2]^*$ is
\begin{equation} \label{pcr}
\phi (\tau_{(12)} \circ R) = (\tau_{(12)} \circ R) \phi,
\end{equation}
which says that $\phi$ commutes with the $R$-permutation action. Since
$H_R$ is presented by the quadratic relation \eqref{hrrele}, $H_R[m]^*$
consists of $\phi$ that satisfy \eqref{pcr} when applying $(\tau_{(12)} \circ R)$ to any components $i,i+1$ for $1 \leq i \leq m-1$.  This proves \eqref{srmhrm}.

Since $S_m$ is completely reducible over $\kk$ ($\kk$ has characteristic zero),
the above means that $V^{\ot m}$ decomposes, as a $S_m$-representation, into
a sum of the form
\begin{equation}
V^{\ot m} \cong \bigoplus_{\lambda} \rho_{\lambda,S_m} \ot V_{m,\lambda,R},
\end{equation}
where $\rho_{\lambda,S_m} \ot V_{m,\lambda,R}$ corresponds to the
$\rho_{\lambda,S_m}$-isotypical part of $V^{\ot m}$, with respect to
the $R$-permutation action of $S_m$, and $V_{m,\lambda,R}$ has trivial
$S_m$-action.  Then, it immediately follows that
$HR_m= \bigoplus \Id \ot \End(V_{m,\lambda,R})$ with respect to this
decomposition, and since the $\rho_{\lambda,S_m}$ are distinct
irreducible representations, \eqref{hrmsrm} follows, and hence also
\eqref{shrm}.  The remaining statements are immediate, and the theorem
is proved.
\end{proof}
\begin{rem}
  The first part of the above theorem, \eqref{srmhrm}, is still
  true if the charactistic of $\kk$ is not zero, or if the condition on
  unitarity is dropped and $S_m$ is replaced by $B_m$, both using the
  same proof as above. However, the next two parts may not
  generalize (since $B_m$ is not finite, its representations are not
  completely reducible in general, so the double commutant arguments
  fail). Nonetheless, the inclusion $\supseteq$ in \eqref{hrmsrm} is still
  true and well-known in all cases (i.e., \eqref{rpae} defines
  endomorphisms of comodules).
\end{rem}
\begin{rem}
  As a special case of the above, we immediately get the usual Schur-Weyl
  duality for $R=1$ (and it is essentially the same as Weyl's original
  proof), except that the usual statement also says that $HR_m$ is
  generated by the diagonal action of $GL(V)$.  To get this last fact,
  as in Weyl's proof, one may use the fact that the symmetric elements
  of $W^{\ot m}$ are generated by the diagonal elements $w^{\ot m}$, for
  any vector space $W$ over an infinite field.
\end{rem}

\footnotesize{
\bibliography{master}
\bibliographystyle{amsalpha}
{\bf T.S.}: Department of
  Mathematics, University of Chicago, 5734 S. University Ave, Chicago IL
  60637, USA;\\
  \hphantom{x}\quad\, {\tt trasched@math.uchicago.edu}}

\end{document}